\input   amstex
\documentstyle{amsppt}
\input bull-ppt
\keyedby{bull267/clh}

\topmatter
\cvol{26}
\cvolyear{1992}
\cmonth{April}
\cyear{1992}
\cvolno{2}
\cpgs{253-257}
\title   Ramanujan   duals   and   automorphic   spectrum  
 \endtitle
\author   M.   Burger,   J.   S.   Li,   and   P.   
Sarnak\endauthor
\shorttitle{Ramanujan   duals}
\address Graduate Center, City University of New York, New 
York 10036\endaddress
\address Department of Mathematics, University of 
Maryland, College Park,
Maryland 20742\endaddress
\address Department of Mathematics, Princeton University,
Princeton, New Jersey 08544 and IBM Research Division, 
Almaden
research center, 650 Harry Rd, San Jos\'e, California 
95120\endaddress
\date April   16,   1991   and,   in   revised   form,   
July 1991\enddate
\abstract We   introduce   the   notion   of   the   
automorphic
dual   of   a   matrix   algebraic   group   defined   
over   $Q$.
This   is   the   part   of   the   unitary   dual   that  
 corresponds
to   arithmetic   spectrum.   Basic   functorial   
properties   of
this   set   are   derived   and   used   both   to   
deduce   arithmetic
vanishing   theorems   of   ``Ramanujan''   type   as   
well   as
to   give   a   new   construction   of   automorphic   
forms.\endabstract
\endtopmatter

\document

Let   $G$   be   a   semisimple   linear   algebraic   
group   defined
over   $\Bbb   Q.$   In   the   arithmetic   theory   of   
automorphic
forms   the   lattice   $\Gamma   =   G(\Bbb   Z)$   and   
its   congrunce
subgroups
$$
\Gamma   (N)   =   \{   \gamma      \in   G(\Bbb   
Z)\colon\   \gamma   \equiv
I(N)\},\quad   \quad      N\in   \Bbb   N$$
play   a   central   role.   A   basic   problem   is   to 
  understand
the   decomposition into irreducibles  of   the   regular  
 representation   of   $G(\Bbb   R)$
on   $L^2   (\Gamma   (N)\backslash      G(\Bbb   R)).$   
In   general this representation  will   not   be   a   
direct   sum
of   irreducibles,   and   for   our   purposes   of   
defining   the
spectrum,   it   is   best   to   use   the   notions   of 
  weak   containment
and   Fell   topology   on   the   unitary   dual   
$\widehat{   G}(\Bbb   R)$ of
the Lie group $G(\Bbb R)$.
(See   \cite{D,   18.1}.) For   any   closed   subgroup
$H$   of   $G   (   \Bbb      R)$   we   define   the   
spectrum   $\sigma      (H
\backslash      G(R))$   to   be   the   subset   of   
$\widehat{   G}
(\Bbb   R)$   consisting   of   all   $\pi   \in   
\widehat{   G}(\Bbb   R)$
that   are   weakly   contained   in   $L^2(   H\backslash 
     G(\Bbb   R)).$
Furthermore,   if   $\widehat{   G}^1   (\Bbb   R)$   is   
the   set   of   
irreducible   spherical   representations,   we   set   
$\sigma   ^1
(H\backslash      G(\Bbb      R))\coloneq   \sigma   
(H\backslash      G(\Bbb   R))\cap
\widehat{   G}^1   (\Bbb   R)$.   When   $H=   \Gamma      
(N)$,   $\sigma   
(\Gamma   (N)\backslash      G(\Bbb   R))$   consists   of 
  all   $\pi   \in
\widehat{   G}   (\Bbb   R)$   occurring   as   
subrepresentations   of
$L^2   (\Gamma      (N)\backslash      G(\Bbb   R))$   as  
 well   as   those
$\pi$'s   that   are   in   wave   packets   of   unitary  
 Eisenstein
series   \cite{La}.   The   latter   occur   only   when   
$\Gamma   \backslash      
G(\Bbb      R)$   is   not   compact.   We   now   
introduce   the   central
object   of   this   note.

\dfn{Definition}   The   automorphic   (resp.   Ramanujan) 
  dual   of
$G$   is   defined   by
$$
\widehat{   G   }   _{\operatorname{Aut}}=\overline{\bigcup
^\infty   _{   N=1}\sigma   (\Gamma   (N)\backslash        
 G(\Bbb   R))},\tag1$$
$$
\widehat{   G   }_{\operatorname{Raman}}=   \widehat{   G  
 }
_{\operatorname{Aut}}\cap   \widehat{   G   }   ^1   (\Bbb 
  R).$$
Here closure is taken in the topological space $\hat 
G(\Bbb R)$.

Thus,   $\widehat{   G   }   _{\operatorname{Aut}}   $   
is   the   smallest
closed   set   containing   all   the   congruence   
spectrum.   Here   is
an   alternative   description   of   $\widehat{   G   }   
_{\operatorname{Raman}}.$
Let   $   G(\Bbb   R)   =   KAN$   be   an   Iwasawa   
decomposition   of   $G(\Bbb   R);$
then   the   theory   of   spherical   functions   
identifies   $\widehat{   G   }^1
(\Bbb   R)$   with   a   subset   of   $\germ   A   _{\Bbb 
  C}^*/W,$   where   $\germ   A   =
\operatorname{Lie}A,$   $W   =   \operatorname{Weyl}(G,   
A).$   Moreover,
the   Fell   topology   on   $\widehat{   G   }^1   (\Bbb  
 R)$   coincides   with
the   topology   of   $\widehat{   G   }^1   (\Bbb   R)$   
viewed   as   a   subset   of
$\germ   A^*_   {   \Bbb   C   }   /W.$   Let   $\Bbb   D$ 
  be   the   ring   of   invariant
differential   operators   on   the   associated   
symmetric   space   $X$.
Then   the   duality   theorem   \cite{GPS}   shows   that 
  the   spectrum
of   $\Bbb   D$   in   $L^2(\Gamma      \backslash      
X),$   say   $\operatorname{Sp}
_\Gamma      (\Bbb   D)   \subset      \germ   A^*_   {   
\Bbb   C   }/W$   is   the   image
of   $\sigma^1(\Gamma   \backslash      G(\Bbb   R))$   in 
  $\germ   A^*
_   {   \Bbb   C   }   /W$   under   the   above   
identification.   In   particular,
$   \widehat{   G   }   _{\operatorname{Raman}}$   is   
identified   with
$$
\overline{\bigcup   ^\infty      _{   N   =   1}   
\operatorname{Sp}
_{\Gamma   (N)}(\Bbb   D)}\subset      \germ   A^*_   {   
\Bbb   C   }/W.$$
That   there   should   be   restrictions   on   
$\widehat{   G   }_   {
\operatorname   {   Raman   }}$   and   $\widehat{   G   
}_{\operatorname{Aut}}$
has   its   roots   in   the   representation   theoretic  
 reinterpretation
of   the   classical   Ramanujan   conjectures   due   to  
 Satake   \cite{Sa}.
Identifying   the   above   sets   may   be   viewed   as  
 the   general   
Ramanujan   conjectures.   For   example,   Selberg's   
1/4-conjecture
may   be   stated   as   follows:   For   $G   =   
\operatorname{SL}_2,$
$$
\widehat{   G   }_   {   \operatorname   {   Raman   }}   
=\{\boldone\}
\cup   \widehat{   G   }^1   (\Bbb   R)   
_{\operatorname{temp}},\tag2$$
where,   in   general,   $\widehat{   G   }(\Bbb   
R)_{\operatorname{temp}}
\coloneq   \sigma   (G(\Bbb   R))$   is   the   set   of   
tempered   representations,
and   $\widehat{   G   }^1   (\Bbb   R)   
_{\operatorname{temp}}=   \widehat{   G   }
(\Bbb   R)_{\operatorname{temp}}   \cap   \widehat{   G   
}^1   (\Bbb   R).$
(See   \cite{CHH}   for   equivalent   definitions   of   
temperedness.)
\enddfn

While   the   individual   sets   $\sigma      (   \Gamma  
 (N)\backslash      G(\Bbb   R))$
are   intractable,   the   set   $\widehat{   G   
}_{\operatorname{Aut}}$
(and   $\widehat{   G   }_   {   \operatorname   {   Raman 
  }})$   enjoy   certain
functorial   properties.

\proclaim{Theorem   1}   Let   $G$   be   a   connected   
semisimple   linear
algebraic   group   defined   over   $\Bbb   Q$   and   $H 
  <   G$   a   $\Bbb   Q$-subgroup

{\rm   (i)}   $\operatorname{Ind}^{G(\Bbb   R)}_{   H(\Bbb 
  R)}   
\widehat{   H}
_{\operatorname{Aut}}\subset      \widehat{   G   
}_{\operatorname{Aut}}.$

{\rm   (ii)}   Assume   that   $H   $   is   semisimple;   
then
$$
\operatorname{Res}_{   H(\Bbb   R)}   \widehat{   G   }_   
{   \operatorname{Aut   }}   
\subset      \widehat{   H}_   {   \operatorname{Aut   }}  
 .$$

{\rm   (iii)}   $\widehat{   G   }   _   {   
\operatorname{Aut   }}   
\otimes      \widehat{   G   }   _   {   \operatorname{Aut 
  }}      \subset   
\widehat{   G   }_   {   \operatorname{Aut   }}   .$
\endproclaim

   A   word   about   the   meaning   of   these   
inclusions.   Firstly,  
Ind
denotes unitary induction
  and   Res   stands   for   restriction.   By
the   inclusion,   say   in   (i),   we   mean   that   if 
  $\pi   '   \in
\widehat{   H}_   {   \operatorname{Aut   }}   $   and   
$\pi$   is   weakly
contained   in   $\operatorname{Ind}^{G(\Bbb   R)}_{   
H(\Bbb   R)}
\pi   '$   then   $\pi   \in   \widehat{   G   }_   {   
\operatorname{Aut   }}   .$
(i)   produces   (after   a   local   calculation)   
elements   in
$\widehat{   G   }   _   {   \operatorname{Aut   }}   $   
from   ones   in
$\widehat{   H}_   {   \operatorname{Aut   }}   $   and   
yields   a   new
method   for   constructing   automorphic   
representations.   
Observe   also   that   if   $\pi   \in   \widehat{   G   
}(\Bbb   R)$
is an  isolated point  then   $\pi   \in   \widehat{   G   
}   _   {   \operatorname{Aut   }}   $
implies   that   $\pi$   occurs   as   a   
subrepresentation   in   $L^2
(\Gamma   (N)\backslash      G(\Bbb   R))$   for   some   
$N$.   This
fact   will   be   used   below   to   construct   certain 
  automorphic
cohomological   representations.   (ii)   allow   one   to 
  transfer
setwise   upper   bounds   on   $\widehat{   H}_   {   
\operatorname{Aut   }}   $
to   $\widehat{   G   }_   {   \operatorname{Aut   }}   $  
 and   for   many   $G$'s
gives   nontrivial   approximations   to   the      
Ramanujan   conjectures.
(iii)   exhibits   a   certain   internal   structure   of 
  the   set
$\widehat{   G   }_   {   \operatorname{Aut   }}   .$   We 
  illustrate
the   use   of   Theorem   1   with   some   examples.

\ex{Example   A}   If   $H   =   \{e\}$   then   (i)   
implies   that
$$
\widehat{   G   }_   {   \operatorname{Aut   }}   
\supset   \widehat{   G(\Bbb   R)}   _{\operatorname{temp}}
\cup   \{\boldone\}.\tag3$$
When   $G(\!\Bbb   Z\!)$   is   cocompact   this   follows 
  also   from
de   George-Wallach   \cite{GW}.   In   comparison   with  
 (2)
one   might   hope   that   (3)   is   an   equality.   
However,
using   other   $H$'s   and   (i)   one   finds   
typically   that
$\widehat{   G   }_   {   \operatorname{Aut   }}   $
contains   nontrivial,   nontempered   spectrum.   For
$G=   \operatorname{Sp}(4)$   the   failure   of   the
naive   Ramanujan   conjecture   has   been   observed   by
Howe   and   Piatetski-Shapiro   \cite{HP-S}   using
theta   liftings.
\endex

\ex{Example   B}   Let   $k   /   \Bbb      Q$   be   a   
totally
real   field,   $q$   a   quadratic   form   over   $k$   
such
that   $q$   has   signature   $(n,   1)$   over   $\Bbb   
R$,
and   all   other   conjugates   are   definite.   Let   $G=
\operatorname{Res}_{k/\Bbb   Q}   \operatorname{SO}(q).$
Then   $G(\Bbb   R)$   is   of   $\Bbb   R$-rank   one   
and   the
noncompact   factor   is   $\operatorname{SO}(n,   1).$
We   identify   $\germ   A^*$   with   $\Bbb   R$   by
sending   $\rho   $   to   $(n   -   1)/   2$.   With   this
normalization   $\widehat{   G   }(\Bbb   R)^1$   is
identified   with   $i\Bbb   R   \cup   [-\rho   ,   \rho  
 ]
\subset      \Bbb   C\operatorname{modulo}\{\pm   1\}.$
\cite{K}.   We   parametrize   $\widehat{   G   }^1   
(\Bbb   R)$
by   $s\in   i   \Bbb   R^+   \cup   [0,   \rho   ]$   and 
  denote
the   corresponding   representation   by   $\pi   _s$.
Let   $\varphi   _0   ,   \dots   ,   \varphi   _n$   be   
an
orthogonal   basis   of   $q   $   such   that   
$q(\varphi   _i)
   >   0,$   $0   \le   i   \le   n   -   1,$   and   $   
q(\varphi   _n)   <   0.$
Define   $H   =\operatorname{Res}   _{k/   \Bbb   Q}\{g\in
\operatorname{SO}(q)\colon\   g(\varphi   _1)=
\varphi      _1\}.$   Applying   Theorem   1(i)   to the 
trivial representation   $\boldone
\in   \widehat{   H}_   {   \operatorname{Aut   }}   $   
we   find   that
$$
\sigma   (H(\Bbb   R)\backslash   G(\Bbb   R))   \subset   
   \widehat{   G   }_   {
\operatorname{Aut   }}   .$$
Now   $\sigma   ^1   (H(\Bbb   R)\backslash   G(\Bbb   R))$
has   been   computed   (\cite{F}),   and   we   find
$$
\widehat{   G   }_   {   \operatorname   {   Raman   }}   
\supset
\{\rho   ,   \rho   -1,   \rho   -2   ,   \dots   \}\cup
i\Bbb   R^+.\tag4$$
In   particular,   for   $n   \ge   4$   there are 
nontrivial nontempered spherical automorphic
representations.

To   find   upper   bounds   on   $\widehat{   G   }_   {  
 \operatorname   {   Raman   }}$
one   uses   Theorem   1(ii)   and   $$H=\operatorname{Res}
_{k/   \Bbb   Q}\{g\in   \operatorname{SO}   (q)\colon
\      g(\varphi   _i)=\varphi   _   i,\ 1   \le   i\le   
n   -4\}.$$
Combining   the   Jacquet-Langlands   correspondence
\cite{JL}   with   the   Gelbart-Jacquet   lift   \cite{GJ}
one   concludes   that
$$
\widehat{   H}_   {   \operatorname   {   Raman   
}}\subset   
i\Bbb   R^+   \cup   [0,\tfrac   12]   \cup   
\{\boldone\}.$$
Applying   (ii)   it   follows    that
$$
\widehat{   G   }_   {   \operatorname   {   Raman   }}   
\subset   
i\Bbb   R^+   \cup   [0,   \rho   -\tfrac12]\cup   \{\rho  
 \}.
\tag   5$$
In   the   special   case   $k   =   \Bbb   Q$,   $n   \ge 
  4$   this
result   has   also   been   obtained   by   \cite{EGM}   
and
\cite{LP-SS}   using   Poincar\'e   series.   Assuming
the   Ramanujan   conjecture   at   $\infty   $   for   
$\operatorname{GL}
(2)$   one   deduces
$$
\widehat{   G   }_   {   \operatorname   {   Raman   }}   
\subset   
i\Bbb   R^+   \cup   [0,   \rho   -1]\cup   \{\rho   
\}.\tag6$$
(Compare   with   (4).)   The   natural   conjecture   
arising
from   (4)   and   (6)   is
$$
\widehat{   G   }_   {   \operatorname   {   Raman   }}   
=   i\Bbb   R
^+\cup   \{\rho   ,   \rho   -1   ,   \dots   \}.$$
This   is   apparently   consistent   with   Arthur's   
conjectures   \cite{A}.
\endex

\ex   {Example   C}   Let   $\Bbb   F_{4(-20)}$   be   the 
  $\Bbb   R$-rank
one   form   of   $\Bbb   F_4.$   Using   a   method   of  
 Borel   \cite{B},
one   may   find   $\Bbb   Q$-groups   $H   <   G$   such  
 that   $G(\Bbb   R),$
$H(\Bbb   R)$   both   have   rank   one,   the   
noncompact   simple
factors   being   $\Bbb   F   _{4(-20)}$   and   
$\operatorname{Spin}(8,   1)$
respectively.   With   notations   similar   to   Example  
 B,   one
may   identify   $\widehat{   G   }^1(\Bbb   R)$   with   
$   i\Bbb   R^+
\cup   [0,   5]   \cup   \{11\},$   here   $\rho   =   
11.$   One   may
compute   $\sigma   ^1   (   H(\Bbb   R)\backslash   
G(\Bbb   R))$
and   using   Theorem   1(i)   find   that
$$
\widehat{   G   }_   {   \operatorname   {   Raman   }}   
\supset
i\Bbb   R^+   \cup   \{3,   11\}.$$
\endex

\ex{Example   D}   Consider   now   $F_{4   (4)},$   the   
split   real
form   of   $F_4.$   The   corresponding   symmetric   
space   has
dimension   28.   For   any   cocompact   lattice   
$\Gamma   \subset   
F_{4(4)}$   one   knows   from   Vogan-Zuckerman   
\cite{VZ}   that
the   Betti   numbers   $\beta   ^m   (\Gamma   )=   0$   
for   $0   <   m   <   8$
or   $20   <   m   <   28, m\neq 4$ or $24$, in these 
latter dimensions
all the cohomology comes from parallel forms of the 
symmetric space.
Nevertheless   using   Theorem   1(i)   we
have   
\endex

\proclaim{Theorem   2}   For   any   cocompact   lattice   
$\Gamma   $   in
$F_{4(4)}$ and $N\ge0$  there   exists   $\Gamma      '   
\subseteq   \Gamma   $   of
finite   index   such   that   $\beta      ^m(\Gamma   
')H\ge N$  for   
$m   =   8,   20.$
\endproclaim
\par
The   proof   of   Theorem   2   makes   use   of   
Matsushima's   formula
\cite{BW}   together   with   a   recent   result   of   
Vogan   ensuring
that   the   unitary   representation   contributing   to  
 the   above
Betti   numbers   is   isolated   in   the   unitary   
dual   of   $F_{4(4)}.$
The   $\Bbb   Q$-subgroup   that   we   use   in   
applying   Theorem   1(i)
   has   real   points   equal   to   
$\operatorname{Spin}(5,   4)$   up
to   compact   factors.   By   the   well-known   result   
of   Oshima-Matsuki
\cite{MO},   we   conclude   that   the   discrete   
series   of   the   
symmetric   space   $F_{4   (4   )}
/\operatorname{Spin}(5,   4)$   contain   a   unitary   
representation
with   nonzero   cohomology   in   degrees   8   and   20, 
  which   is   isolated
in   the   unitary   dual.   The   fact   that   we   have 
  dealt   with   every
lattice   in   $F   _{   4   (   4   )}$   follows   from  
 Margulis's
arithmeticity   theorem   \cite{M},   together   with   the
classification   of   algebraic   groups   over   number   
fields
\cite{T}.

This   method   of   constructing   cohomology   is   
rather   general.
If   $\pi$   is   isolated   in   $\widehat{   G   }(\Bbb  
 R)$   and   is
contained   in   the   automorphic   dual   of   $G$   
then   it   occurs
discretely   in   $L^2(\Gamma      \backslash      G(\Bbb  
 R))$   for
$\Gamma   $   a   congruence   subgroup   of   deep   
enough   level.
David   Vogan   has   recently   obtained   the   
necessary   and
sufficient   conditions   for   a   unitary   representation
with   nonzero   cohomology   to   be   isolated,   which  
 implies
that   most   of   them   do.   Theorem   1   then   
allows   us   to   obtain
nonvanishing   of   cohomology   in   a   large   number   
of   cases.

To   end,   we   remark   that   these   ideas   extend   
in   a   natural
way   to   $S$-arithmetic   groups.   The   proof   of   
Theorem   1(i)
consists   of   approximating,   in   a   suitable   way,  
 congruence
subgroups   of   $H(\Bbb   Z)$   by   congruence   
subgroups   of
$G(\Bbb   Z)$   and   then   applying   criteria   of   
weak   containment.
For   the   proofs   of   Theorem   1(ii),   (iii)   we   
refer   the
reader   to   \cite{BS}.

\heading   Acknowledgments\endheading   We   would   like  
 to   thank
David   Vogan   for   sharing   his   insights   into   
unitary
representations   and   F.   Bien   for   interesting   
conversations.

\enddocument